# DISCUSSION OF "EQUI-ENERGY SAMPLER" BY KOU, ZHOU AND WONG

By Ming-Hui Chen and Sungduk Kim

*University of Connecticut*

**1. Introduction.** We first would like to congratulate the authors for their interesting paper on the development of the innovative equi-energy (EE) sampler. The EE sampler provides a solution, which may be better than existing methods, to a challenging MCMC sampling problem, that is, sampling from a multimodal target distribution $\pi(x)$. The EE sampler can be understood as follows. In the equi-energy jump step, (i) points may move within the same mode; or (ii) points may move between two modes; but (iii) points cannot move from one energy ring to another energy ring. In the Metropolis–Hastings (MH) step, points move locally. Although in the MH step, points may not be able to move freely from one mode to another mode, the MH step does help a point to move from one energy ring to another energy ring locally. To maintain certain balance between these two types of operations, an EE jump probability $p_{ee}$ must be specified. Thus, the MH move and the equi-energy jump play distinct roles in the EE sampler. This unique feature makes the EE sampler quite attractive in sampling from a multimodal target distribution.

**2. Tuning and "black-box."** The performance of the EE sampler depends on the number of energy and temperature levels, $K$, energy levels $H_0 < H_1 < \cdots < H_K < H_{K+1} = \infty$, temperature ladders $1 = T_0 < T_1 < \cdots < T_k$, the MH proposal distribution, the proposal distribution used in the equi-energy jump step and the equi-energy jump probability $p_{ee}$. Based on our experience in testing the EE sampler, we felt that the choice of the $H_k$, the MH proposal and $p_{ee}$ are most crucial for obtaining an efficient EE sampler. In addition, the choice of these parameters is problem-dependent. To achieve fast convergence and good mixing, the EE sampler requires extensive tuning of $K$, $H_k$, MH proposal and $p_{ee}$ in particular. A general sampler is designed to be "black box" in the sense that the user need not tune the sampler to the problem. Some attempts have been made for developing such "black-box" samplers in the literature. Neal [4] developed variations on slice sampling that can be used to sample from any continuous distributions and that require little or no tuning. Chen and Schmeiser [2] proposed the random-direction interior-point (RDIP) sampler. RDIP samples from the uniform distribution defined over the region $U = \{(x, y) : 0 < y < \pi(x)\}$ below the curve of the surface defined by $\pi(x)$, which is essentially the same idea used in slice sampling.

**3. Boundedness.** It is not clear why the target distribution $\pi(x)$ must be bounded. Is this a necessary condition required in Theorem 2? It appears that the condition $\sup_x \pi(x) < \infty$ is used only in the construction of energy levels $H_k$ for $k > 0$ for convenience. Would it be possible to relax such an assumption? Otherwise, the EE sampler cannot be applied to sampling from an unbounded $\pi(x)$ such as a gamma distribution with shape parameter less than 1.

If we rewrite
$$D_j = \{x : h(x) \in [H_j, H_{j+1})\} = \{x : \pi(x) \in (\exp(-H_{j+1}), \exp(-H_j)]\},$$
we can see that $D_0$ corresponds to the highest-density region. Thus, if $H_1$ is appropriately specified, and the guideline given in Section 3.3 is applied to the choice of the rest of the $H_j$'s, the boundedness assumption on $\pi(x)$ may not be necessary.

**4. Efficiency.** The proposed EE sampler requires $K(B + N)$ iterations before it starts the lowest-order chain $\{X_n^{(0)}, n \geq 0\}$. Note that here $B$ is the number of "burn-in" iterations and $N$ is the number of iterations used in constructing an empirical energy ring $\hat{D}_j^k$. As it is difficult to determine how quickly a Markov chain $\{X_n^{(k)}\}$ converges, a relatively large $B$ may be needed. If the chain $X^{(k)}$ does not converge, the acceptance probability given in Section 3.1 for the equi-energy move at energy levels lower than $k$ may be problematic. Therefore, the EE sampler is quite inefficient as a large number of "burn-in" iterations will be wasted. This may be particularly a problem when $K$ is large. Interestingly, the authors never disclosed what $B$ and $N$ were used in their illustrative examples. Thus, the choice of $B$ and $N$ should be discussed in Section 3.3.

**5. Applicability in high-dimensional problems.** Based on the guideline of the practical implementation provided in the paper, the number of energy levels $K$ could be roughly proportional to the dimensionality of the target distribution. Thus, for a high-dimensional problem, $K$ could be very large. As a result, the EE sampler may become more inefficient as more "burn-in" iterations are required and at the same time, it may be difficult to tune the parameters involved in the EE sampler.

For example, consider a skewed link model for binary response data proposed by Chen, Dey and Shao [1]. Let $(y_1, y_2, \ldots, y_n)'$ denote an $n \times 1$ vector of $n$ independent dichotomous random variables. Let $x_i = (x_{i1}, \ldots, x_{ip})'$ be a $p \times 1$ vector of covariates. Also let $(w_1, w_2, \ldots, w_n)'$ be a vector of independent latent variables. Then, the skewed link model is formulated as follows: $y_i = 0$ if $w_i < 0$ and 1 if $w_i \geq 0$, where $w_i = x_i'\beta + \delta z_i + \varepsilon_i$, $z_i \sim G$, $\varepsilon_i \sim F$, $z_i$ and $\varepsilon_i$ are independent, $\beta = (\beta_1, \ldots, \beta_p)'$ is a $p \times 1$ vector of regression coefficients, $\delta$ is the skewness parameter, $G$ is a known cumulative distribution function (c.d.f.) of a skewed distribution, and $F$ is a known c.d.f. of a symmetric distribution. To carry out Bayesian inference for this binary regression model with a skewed link, we need to sample from the joint posterior distribution of $((w_i, z_i), i = 1, \ldots, n, \beta, \delta)$ given the observed data $D$. The dimension of the target distribution is $2n + p + 1$. When the sample size $n$ is large, we face a high-dimensional problem. Notice that the dimension of the target distribution can be reduced considerably if we integrate out $(w_i, z_i)$ from the likelihood function. However, in this case, the resulting posterior distribution $\pi(\beta, \delta | D)$ contains many analytically intractable integrals, which could make the EE sampler expensive or even infeasible to implement. The skewed link model is only a simple illustration of a high-dimensional problem. Sampling from the posterior distribution under nonlinear mixed-effects models with missing covariates considered in [5] could be even more challenging.

In contrast, the popular Gibbs sampler may be more attractive and perhaps more suitable for a high-dimensional problem because the Gibbs sampler requires only sampling from low-dimensional conditional distributions. As MH sampling can be embedded into a Gibbs step, would it be possible to develop an EE-within Gibbs sampler?

**6. Statistical estimation.** In the paper, the authors proposed a sophisticated but interesting Monte Carlo method to estimate the expectation $E_{\pi_0}[g(X)]$ under the target distribution $\pi_0(x) = \pi(x)$ using all chains from the EE sampler. Due to the nature of the EE sampler, the state space $\mathcal{X}$ is partitioned according to the energy levels, that is, $\mathcal{X} = \bigcup_{j=0}^{K} D_j$. Thus, this may be an ideal scenario for applying the partition-weighted Monte Carlo method proposed by Chen and Shao [3]. Let $\{X_i^{(0)}, i = 1, 2, \ldots, n\}$ denote the sample under the chain $X^{(0)}$ ($T = 1$). Then, the partition-weighted Monte Carlo estimator is given by

$$\hat{E}_{\pi_0}[g(X)] = \frac{1}{n} \sum_{i=1}^{n} \sum_{j=0}^{K} w_j g(X_i^{(0)}) 1\{X_i^{(0)} \in D_j\},$$

where the indicator function $1\{X_i^{(0)} \in D_j\} = 1$ if $X_i^{(0)} \in D_j$ and 0 otherwise, and $w_j$ is the weight assigned to the $j$th partition. The weights $w_j$ may be estimated using the combined sample, $\{X^{(k)}, k = 1, 2, \ldots, K\}$, under the $\pi_k$ for $k = 1, 2, \ldots, K$.

**7. Example 1.** We consider sampling from a two-dimensional normal mixture,

$$(7.1) \qquad f(x) = \sum_{i=1}^{2} \frac{1}{2} \left[ \frac{1}{2\pi} |\Sigma_i|^{-1/2} \exp\left\{ -\frac{1}{2}(x - \mu_i)' \Sigma_i^{-1}(x - \mu_i) \right\} \right],$$

where

$$x = (x_1, x_2)', \qquad \mu_1' = (0, 0), \qquad \mu_2' = (5, 5)$$

and

$$\Sigma_i = \begin{pmatrix} \sigma_1^2 & \sigma_1 \sigma_2 \rho_i \\ \sigma_1 \sigma_2 \rho_i & \sigma_2^2 \end{pmatrix}$$

with $\sigma_1 = \sigma_2 = 1.0$, $\rho_1 = 0.99$ and $\rho_2 = -0.99$. The purpose of this example is to examine performance of the EE sampler under a bivariate normal distribution with a high correlation between $X_1$ and $X_2$. Since the minimum value of the energy function $h(x) = -\log(f(x))$ is around $\log(4\pi\sigma_1\sigma_2\sqrt{1.0 - \rho_i^2}) \approx 0.573$, we took $H_0 = 0.5$. $K$ was set to 2. The energy ladder was set between $H_{\min}$ and $H_{\min} + 100$ in a geometric progression, and the temperatures were between 1 and 60. The equi-energy jump probability $p_{\text{ee}}$ was taken to be 0.1. The initial states of the chain $X^{(i)}$ were drawn uniformly from $[0, 1]^2$. The MH proposal was taken to be bivariate Gaussian: $X_{n+1}^{(i)} \sim N_2(X_n^{(i)}, \tau_i^2 T_i I_2)$, where the MH proposal step size $\tau_i$ for the $i$th-order chain $X^{(i)}$ was taken to be 0.5 such that the acceptance ratio was in the range of $(0.23, 0.29)$. The overall acceptance rate for the MH move in the EE sampler was 0.26. We used 2000 iterations to burn in the EE sampler and then generated 20,000 iterations. Figure 1 shows autocorrelations and the samples generated in each chain based on the last 10,000 iterations. We can see, from Figure 1, that the EE sampler works remarkably well and the high correlations do not impose any difficulty for the EE sampler at all.

**8. Example 2.** In this example, we consider another extreme and more challenging case, in which we assume a normal mixture distribution with different variances. Specifically, in (7.1) we take

$$\Sigma_i = \begin{pmatrix} \sigma_{i1}^2 & \sigma_{i1}\sigma_{i2}\rho_i \\ \sigma_{i1}\sigma_{i2}\rho_i & \sigma_{i2}^2 \end{pmatrix}$$

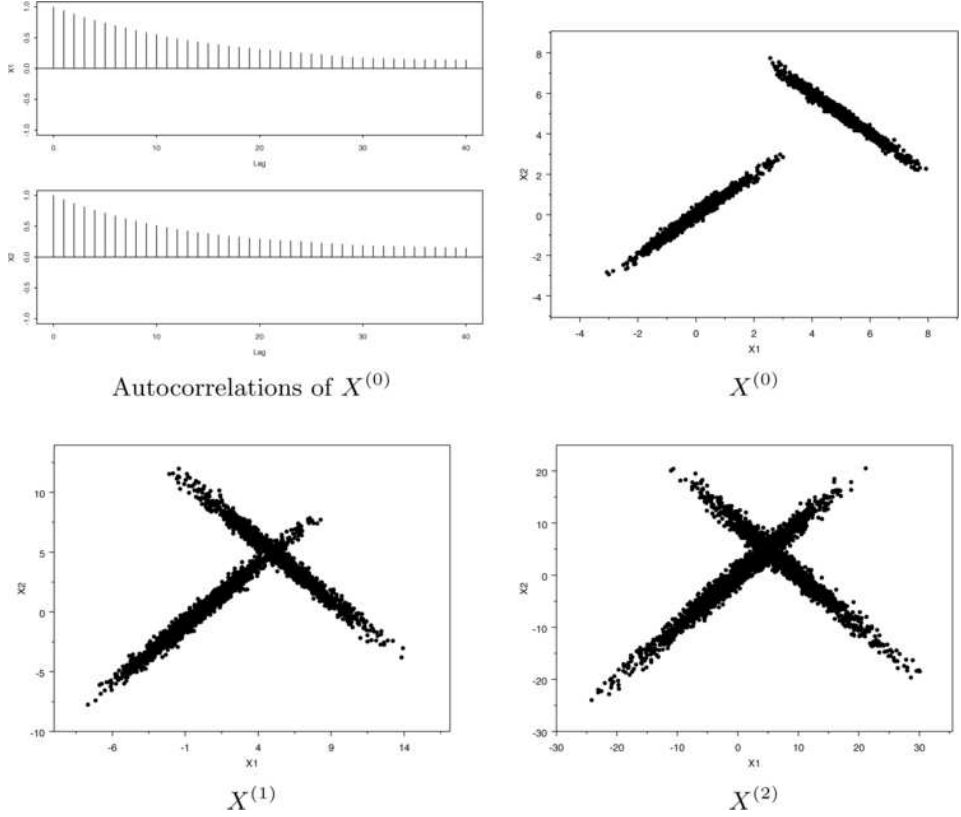

Fig. 1. *Plots of EE samples from a normal mixture distribution with equal variances.*

with $\sigma_{11} = \sigma_{12} = 0.01$, $\sigma_{21} = \sigma_{22} = 1.0$ and $\rho_1 = \rho_2 = 0$. Since the minimum value of the energy function $h(x)$ is around $-6.679$, we took $H_0 = -7.0$. We first tried the same setting for the energy and temperature ladders with $K = 2$, $p_{\text{ee}} = 0.1$ and the MH proposal $N_2(X_n^{(i)}, \tau_i^2 T_i I_2)$. The chain $X^{(0)}$ was trapped around one mode and did not move from one mode to another at all. A similar result was obtained when we set $K = 4$. So, it did not help to simply increase $K$. One potential reason for this may be the choice of the MH proposal $N_2(X_n^{(0)}, \tau_0^2 I_2)$ at the lowest energy level. If $\tau_0$ is large, a candidate point around the mode with a smaller variance is likely to be rejected. On the other hand, the chain with a small $\tau_0$ may move more frequently, but the resulting samples will be highly correlated.

Intuitively, an improvement could be made by increasing $K$, tuning energy and temperature ladders, choosing a better MH proposal and a more appropriate $p_{\text{ee}}$. Several attempts along these lines were made to improve the EE sampler and the results based on one of those trials are given below. In this attempt, $K$ was set to 6, and $H_1 = \log(4\pi) + \alpha = 2.53 + \alpha$, where $\alpha$ was set to 0.6. The energy ladder was set between $H_1$ and $H_{\min} + 100$ in a geometric progression, the temperatures were between 1 and 70, and $p_{ee} = 0.5$. The MH proposals were specified as $N_2(X_n^{(i)}, \tau_i^2 T_i I_2)$ for $i > 0$ and $N_2(\mu(X_n^{(0)}), \Sigma(X_n^{(0)}))$ at the lowest energy level, where $\mu(X_n^{(0)})$ was chosen to be the mode of the target distribution based upon the location of the current point $X_n^{(0)}$ and $\Sigma(X_n^{(0)})$ was specified in a similar fashion as $\mu(X_n^{(0)})$. We used 20,000 iterations to burn in the EE sampler and then generated 50,000 iterations. Figure 2 shows the plots of the samples generated in $X^{(0)}$ based on all 50,000 iterations. The resulting chain had excellent mixing around each mode, and the chain also did move from one mode to another mode. However, the chain did not move as freely as expected.

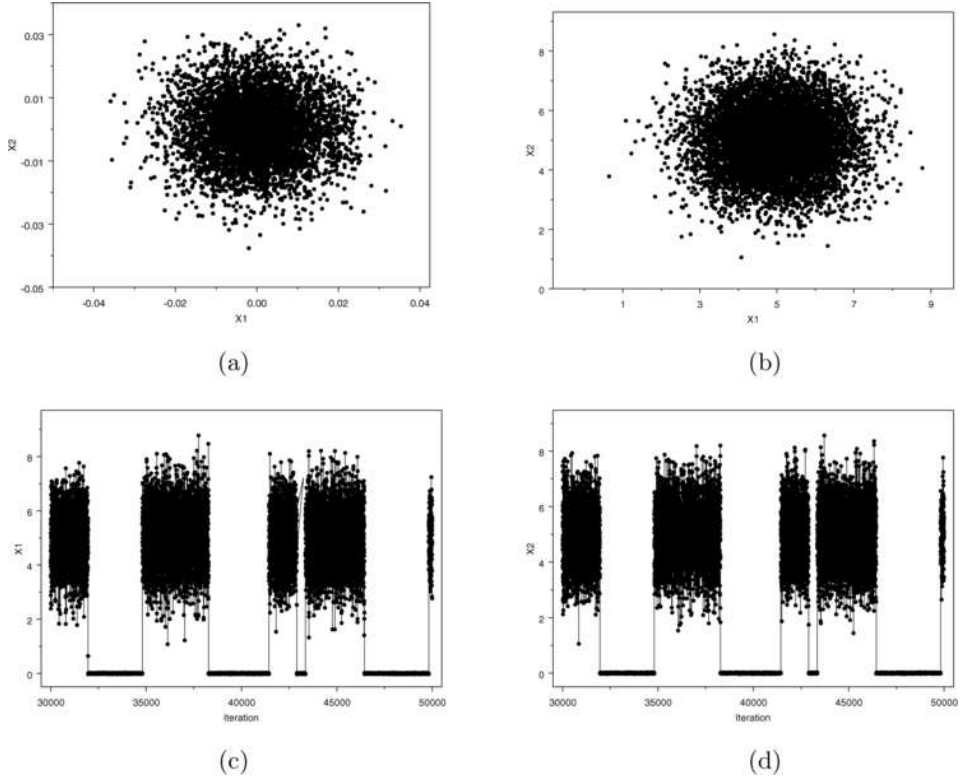

FIG. 2. *Normal mixture distribution with unequal variances. Samples of* $X^{(0)} = (X_1^{(0)}, X_2^{(0)})$ *around* (a) *mode* (0,0) *and* (b) *mode* (5,5). *The marginal sample paths of* $X_1^{(0)}$ (c) *and* $X_2^{(0)}$ (d).

Due to lack of experience in using the EE sampler, we are not sure at this moment whether the EE sampler can be further improved for this example. If so, we do not know how. We would like the authors to shed light on this.

**9. Discussion.** The EE sampler is a potentially useful and effective tool for sampling from a multimodal distribution. However, as shown in Example 2, the EE sampler did experience some difficulty in sampling from a bivariate normal distribution with different variances. For the unequal variance case, the guidelines for practical implementation provided in the paper may not be sufficient. The statement, "the sampler can jump freely between the states with similar energy levels," may not be accurate as well.

As a uniform proposal was suggested for the equi-energy move, it becomes apparent that the points around the modes corresponding to larger variances are more likely to be selected than those corresponding to smaller variances. Initially, we thought that an improvement might be made by assigning a larger probability to the points from the mixand with a smaller variance. However, this would not work as the resulting acceptance probability would become small. Thus, a more likely selected point may be less likely to be accepted. It does appear that a uniform proposal may be a good choice for the equi-energy move.

## REFERENCES

[1] CHEN, M.-H., DEY, D. K. and SHAO, Q.-M. (1999). A new skewed link model for dichotomous quantal response data. *J. Amer. Statist. Assoc.* **94** 1172–1186. MR1731481
[2] CHEN, M.-H. and SCHMEISER, B. W. (1998). Toward black-box sampling: A random-direction interior-point Markov chain approach. *J. Comput. Graph. Statist.* **7** 1–22. MR1628255
[3] CHEN, M.-H. and SHAO, Q.-M. (2002). Partition-weighted Monte Carlo estimation. *Ann. Inst. Statist. Math.* **54** 338–354. MR1910177
[4] NEAL, R. M. (2003). Slice sampling (with discussion). *Ann. Statist.* **31** 705–767. MR1994729
[5] WU, L. (2004). Exact and approximate inferences for nonlinear mixed-effects models with missing covariates. *J. Amer. Statist. Assoc.* **99** 700–709. MR2090904

DEPARTMENT OF STATISTICS
UNIVERSITY OF CONNECTICUT
STORRS, CONNECTICUT 06269-4120
USA
E-MAIL: mhchen@stat.uconn.edu
        sdkim@stat.uconn.edu